\documentclass[11pt]{amsart}
\usepackage{latexsym,amssymb,amsmath}
\usepackage{epsfig}

\newtheorem{lemma}{Lemma}
\newtheorem{proposition}{Proposition}
\newtheorem{theorem}{Theorem}

\theoremstyle{definition}
\newtheorem{definition}{Definition}
\newtheorem{example}{Example}
\newtheorem{remark}{Remark}
\makeindex

\def\impart{{\sf Im}}
\def\realpart{{\sf Re}}

\def\s{\sigma}

\def\tcos{\rm cos}\def\tsin{\rm sin}

\def\frp#1{\frac{\partial}{\partial{#1}}}
\def\frpp#1#2{\frac{\partial{#1}}{\partial{#2}}}

\def\xx#1{x^{#1}}

\def\BC{\mathbb C}\def\BO{\mathbb O}
\def\BR{\mathbb R}\def\BE{\mathbb E}

\def\ta#1{{\theta}^{#1}}
\def\tdim{\rm dim}
\def\hd{,...,}
\def\ww{\wedge}
\def\upperp{{}^\perp}

\def\ppp#1#2{p^{#1}_{#2}}

\def\cI{{\mathcal I}}

\def\ci{{\mathcal I}}
\def\cV{{\mathcal V}}

\def\11{\mathbf 1}

\def\FF{\mathbb F}

\def\l{\lambda}
\def\a{\alpha}

\def\b{\beta}

\def\s{\sigma}

\def\th{\theta}

\def\ra{{\mathord{\;\rightarrow\;}}}

\def\xx#1{x^{#1}}
\def\uu#1{u^{#1}}

\def\tdim{\text{dim}\,}\def\tvol{\text{vol}\,}
\def\tcodim{\text{codim}\,}

\def\tker{\text{ker}\,}

\def\tmod{\text{ mod }}

\begin{document}

\title{ EXTERIOR DIFFERENTIAL SYSTEMS AND BILLIARDS} 
\author{J.M. Landsberg}
\address{Department of Mathematics,
Texas A\& M University,
Mailstop 3368,
College Station, TX 77843-3368, USA}

\begin{abstract}
I describe work in progress with Baryshnikov
and Zharnitsky on periodic billiard orbits that leads one
to an  exterior differential system (EDS). I then give a
brief introduction to EDS illustrated by several examples. 
\end{abstract}

\maketitle

\section{Introduction}

The purpose of these notes is to introduce the reader to the techniques
of exterior differential systems (EDS) in the context of a problem in billiards.
The approach in this article is different from that of \cite{IvL} and
\cite{Ldaewoo}, which begin  with a study of
linear Pfaffian systems, an important special class of EDS.
The billiard problem results in an EDS that is not a linear Pfaffian
system, so   these notes  deal immediately general EDS.
For the interested reader,   two   references regarding
EDS are \cite{BCG3} and \cite{IvL}. The first is a definitive reference
and the second contains an introduction to the subject via
differential  geometry.
For more details about anything regarding EDS the reader
can consult either of these two sources.
   Cartan's book on EDS  \cite{C1945}  is still worth looking at,
especially   the second half, which is a series of beautiful
examples.

We generally will work in the real analytic category, although
all the non-existence results discussed here imply non-existence
of smooth solutions.  

\subsection*{Notation} If $M$ is a differentiable manifold we
let $TM,T^*M$ denote its tangent and cotangent bundles,
$\Omega^d(M)$ the set of differential forms on $M$ of degree
$d$ and $\Omega^*(M)=\oplus_d\Omega^d(M)$. If $I\subset T^*M$
is a subbundle (more precisely, subsheaf), then we
let $\{I\}_{diff}\subset \Omega^*(M)$ denote the differential
ideal generated by $I$, i.e, all elements of $\Omega^*(M)$ of
the form $\a\ww\phi +d\b\ww\psi$ where $\a,\b\in I$ and $\phi,\psi\in\Omega^*(M)$.
$\{ v_1\hd v_n\}$ denotes the linear span of the vectors
$v_i$ if they are vectors, and the subbundle of $\Omega^1(M)$
they generate
if they are one-forms.

\section{Origin of the billiard problem}
Let $D\subset \BR^2$ be a convex domain with its flat
metric. Let $\Delta$ denote the standard Laplacian on
$D$. Then Weyl \cite{weyl} conjectured and Ivrii \cite{ivrii}
proved
$$
\{
{\rm number\ of\ eigenvalues\ of \ }\Delta \  
\leq \l^2
\}
=\frac 1{\pi}{\rm area}(D)\l^2 \pm\frac 14
{\rm length}(\partial D)\l
+o(\l)
$$
where more precisely, Weyl proved the first term
is indeed the leading term and Ivrii proved
the correction term ($+$ with Dirchlet, $-$ with Neumann
boundary conditions), but subject to the following possibly
extraneous hypothesis:

\smallskip

{\it That there does not exist a two parameter family
of periodic billiard trajectories in $D$.
} 

\smallskip

In fact Weyl and Ivrii work in $n$ dimensions but we
have restricted to $n=2$ for notational simplicity.
Also, Ivrii's actual restriction was that there was
not a set of positive measure of periodic billiard
trajectories in the space of all trajectories, but
for the problem at hand, that is equivalent to
the statement above, as remarked in \cite{rychlik}.

I will report on joint work with Y. Barishnikov and
V. Zharnitsky investigating whether this additional
hypothesis is actually necessary or not. But first,
I must explain the hypothesis.

\section{Billiards}
Let $C\subset \BR^2$ be a smooth   curve.
A {\it billiard trajectory} is defined by a particle traveling in
straight lines in the interior of $C$ and  reflecting
at the boundary subject to the law that the angle of incidence
with the tangent line to the curve equals the angle of 
reflection.

 \begin{figure}[!htb]\begin{center}
 \includegraphics[scale=.1]{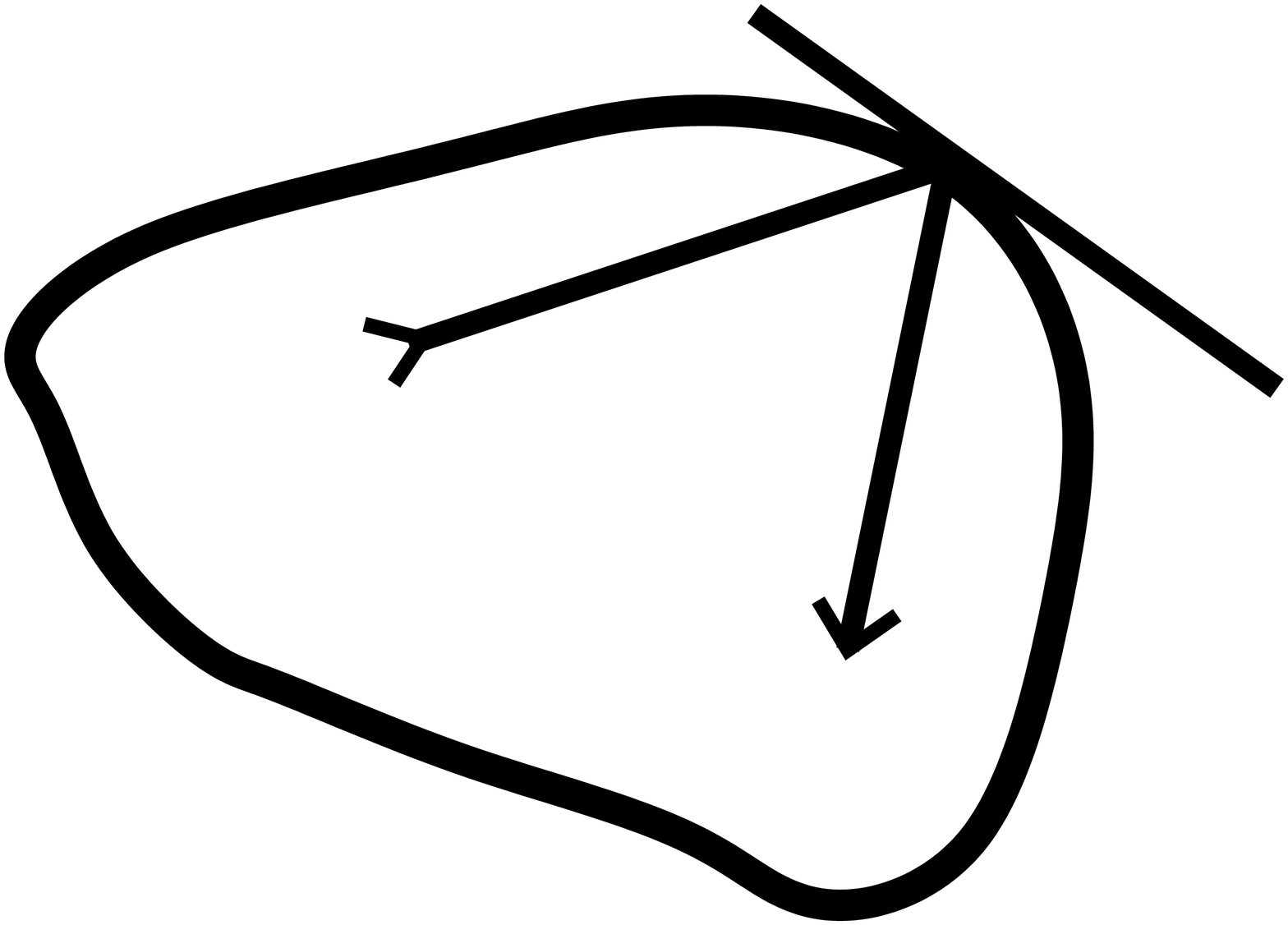}
 \caption{\small{}}  
 \end{center}
 \end{figure}

A trajectory is {\it periodic} if it closes up and
repeats itself. The number of collisions it has with the
boundary of $C$  before repeating is called its
{\it period}.

For example, if $C$ is a circle, then there are many
periodic trajectories.
\begin{figure}[!htb]\begin{center}
\includegraphics[scale=.13]{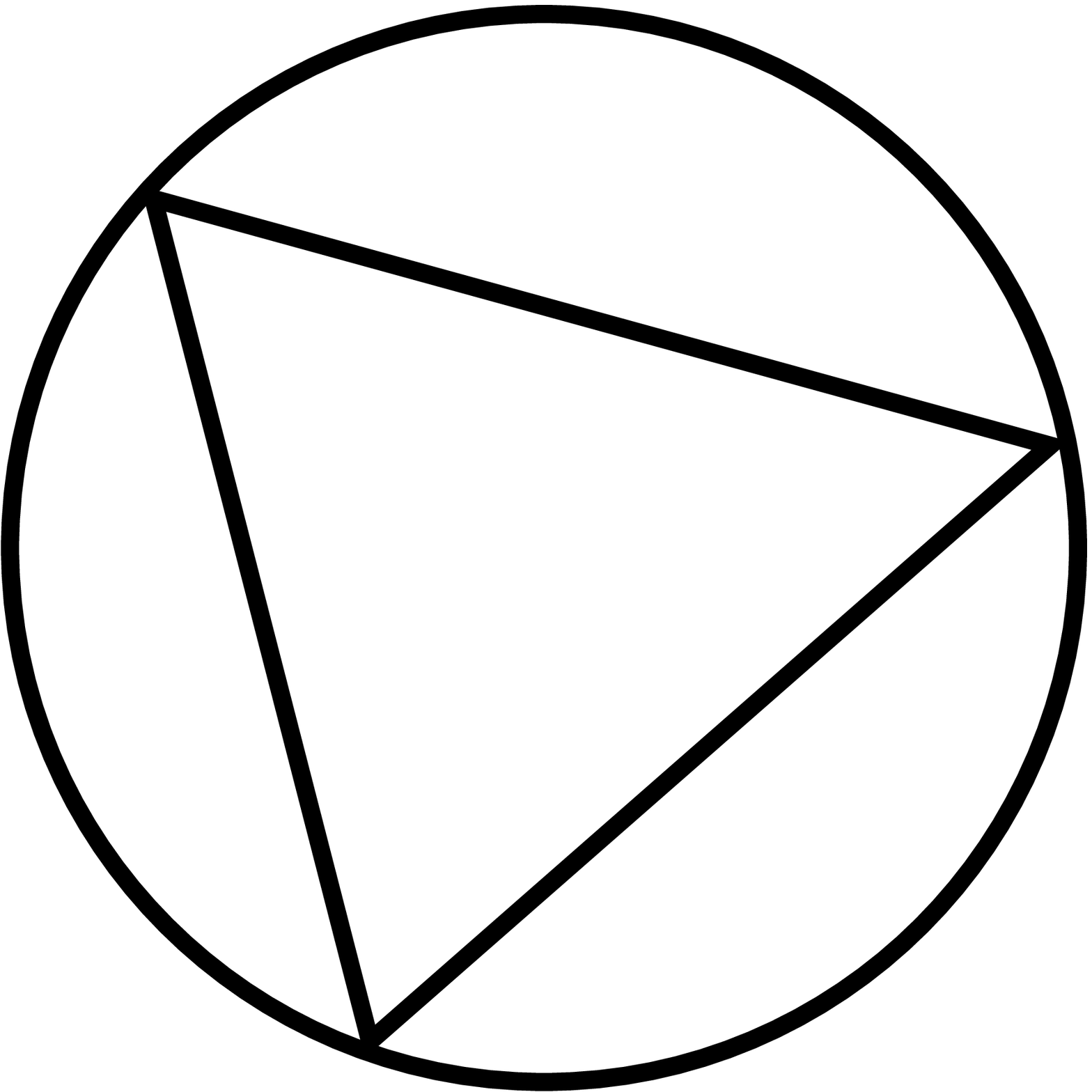}
\caption{\small{}}  
\end{center}
\end{figure}

 Moreover, given a periodic
trajectory in the circle
one can construct a one parameter family of such by
varying the initial point and keeping the angle
constant. It is also true that given an ellipse
and a periodic trajectory on it, one can still obtain
a one parameter family of periodic trajectories 
if one moves the angle just right as one displaces
the initial point. 

\begin{figure}[!htb]\begin{center}
\includegraphics[scale=.15]{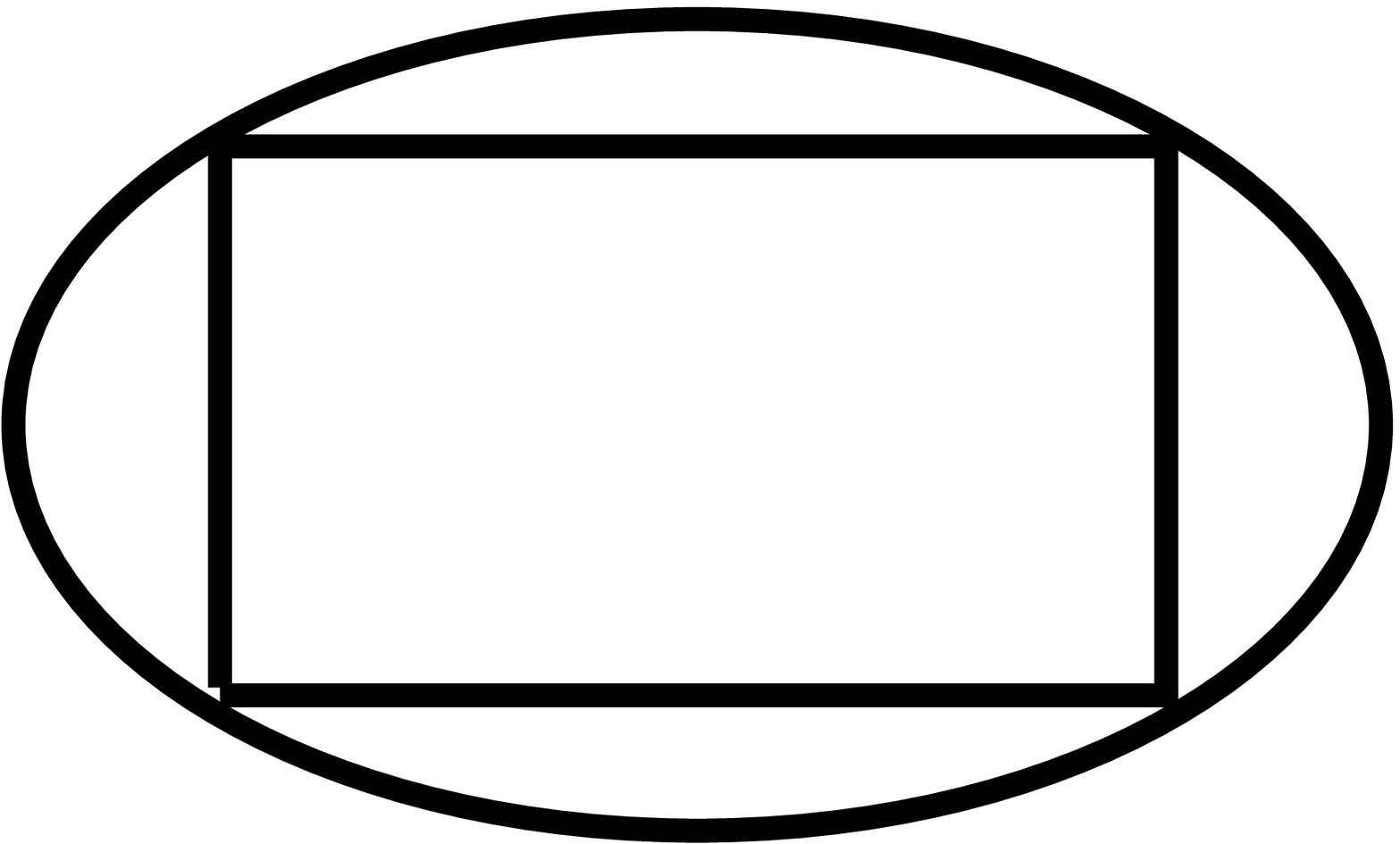}
\caption{\small{}}  
\end{center}
\end{figure}

One can locally parametrize the space of trajectories
 by putting a local parameter on the curve
(e.g. arclength) and measuring the angle of the
trajectory with the tangent line to the curve
(one thinks
of shooting out a trajectory from that point).
In particular, the set of trajectories is a two dimensional space
  and
the existence of a two parameter family of 
periodic trajectories would mean
that there is some point on the curve $C$ such that
no matter what small perturbation of the initial angle
and initial point one makes, the resulting trajectory
is still periodic.

\smallskip

{\it Sound  preposterous?}

\smallskip

Ivrii thought so. In fact, 
legend has it that Ivrii was attempting to prove the
correction term to Weyl's formula and realized it
he could prove it under the assumption
that there are no   periodic billiard trajectories in the domain.
  Fortunately for him (he thought), he was at Moscow State University,
where there were many world experts on billiards.
Allegedly he went in to ask them if
there could be such a curve - they quickly answered:
\lq\lq Of course not!\rq\rq , so he said \lq\lq Great! may I
please have a proof?\rq\rq - they said certainly.
They had  trouble coming up with a full argument immediately
so they told him to come back later in the afternoon. He
returned later that afternoon and they told him that
perhaps it would be better to return the next day...
then it became the next week, ... All this
was nearly 30 years ago and the question is pretty much
as open today as it was then.

\smallskip

Some things are known: for a similar problem  the
answer is that there {\it are} such things: there exist
compact surfaces with Riemannian metrics, all of
whose geodesics are closed \cite{guill}. These 
are called {\it Zoll surfaces} and there are
more of them than was originally expected.

The progress on Ivrii's question is as follows: we
may break it up into a series of questions based on
the period of the trajectory. It is easy to
see that there can be at most a $1$-parameter
family of two-periodic trajectories. (Hint:
what happens when you change the angle a little?) 

In 1989 Rychlik \cite{rychlik} proved that
there are no curves supporting an
open set of $3$-periodic trajectories.
Now there are three published proofs
of Rychlik's theorem \cite{rychlik,stojanov,mpw}, and in these lectures
I will give you a fourth. The four periodic
case is still open and that is the subject 
of my current research with Baryshnikov and
Zharnitsky.

\section{Setting up the problem}\label{setsect}
(The results in this section are joint work
with Baryshnikov and Zharnitsky.)
The problem is local. If we want an $n$-periodic
trajectory, we only need $n$ bits of curve. We can
later close up the bits any way we please (as long
as it closes convexly).

Let $z_1\hd z_n\in \BR^2$. We want to
construct $n$ (germs of) curves, one passing through each point.
The  initial points determine an initial $n$-gon which in turn
tells us what the tangent lines to the curves 
must be at the $z_i$. I.e., the $n$ points
immediately determine the zero-th and
first order terms of the Taylor series for the curves.

Let
$$
N_i=\frac{z_i-z_{i+1}}{|z_i-z_{i+1} |} -
\frac{z_i-z_{i-1}}{|z_i-z_{i-1}|}
$$
and note that $N_i$ points in the direction of
the tangent line to the curve   we are trying to construct.
Let $n_i=N_i/|N_i|$. Let $Jn_i$ denote
the rotation counterclockwise of $n_i$ by $\pi/2$.
We have the following picture

\begin{figure}[!htb]\begin{center}
\includegraphics[scale=1]{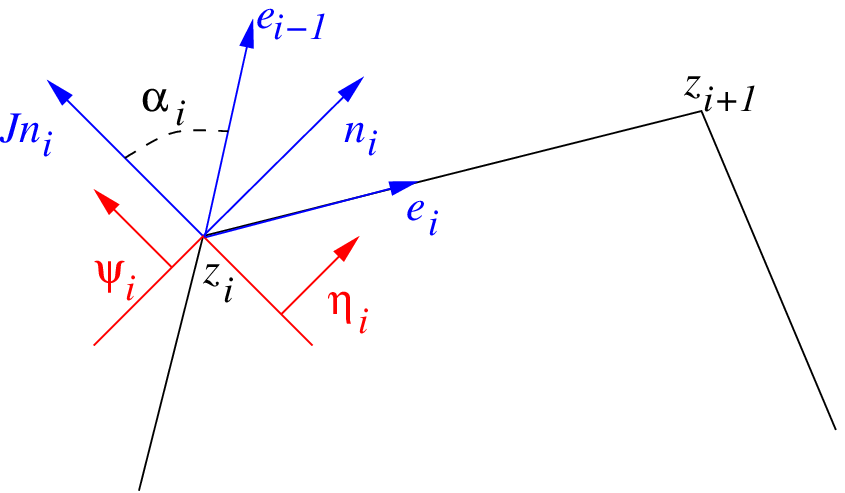}
\caption{\small{}} \label{notation}
\end{center}
\end{figure}

The tangent line 
at $z_i$ must be perpendicular to $J{n_i}$.
Let $\Sigma=(\BR^2)^{\times n}$ denote
the na\"\i ve  configuration space 
(the actual space is an open  subset of this)
where $p=(z_1\hd z_n)\in \Sigma$ is our initial
point.
Define
$$
\psi_i:=\langle Jn_i, dz_i\rangle\in \Omega^1(\Sigma)
$$
and for future reference set $\eta^i=\langle n_i,dz_i\rangle$,
let $\a_i$   be the angle
between $e_{i-1}$ and $Jn_i$, and let $l_i$ be the length of the section
from $z_i$ to $z_{i+1}$.

We have a {\it distribution} $\Delta$ on
$\Sigma$, namely 
$$\Delta=\tker\{\psi_1\hd\psi_n\}\subset T\Sigma
$$

Any two parameter family of $n$-periodic 
trajectories corresponds
to an immersed surface $M^2\hookrightarrow \Sigma$
which is
everywhere tangent to $\Delta$  and subject to
some additional genericity conditions. More precisely, we have

\begin{proposition}\cite{BLZ}\label{blzprop} There exists a one to one correspondence
between (segments of) curves admitting an open subset of
  $n$-periodic trajectories
and immersed
surfaces $i:M^2\ra \Sigma$ tangent to $\Delta$ satisfying
\begin{enumerate}
\item no two points coincide
\item no three points are colinear
\item Any two consecutive points \lq\lq move
independently\rq\rq\  in the manner made
precise in condition \eqref{indepcond} 
described below.
\end{enumerate}
\end{proposition}

Note that the first two conditions are zero-th
order conditions regarding the initial point in
$\Sigma$ about which we want to construct the
surface. The last is a first order condition
which may be described as follows:

  Note
that $(\eta^1_p\hd \eta^n_p,\psi^1_p\hd \psi^n_p)$ gives a basis of $T^*_p\Sigma$
and that this basis varies smoothly - one says
$(\eta^i,\psi^i)$ form a {\it coframing} of 
  of $\Sigma$.
The precise form of condition (3) for   proposition
\ref{blzprop} is that
\begin{equation}\label{indepcond}
i^*(\eta^i\ww\eta^{i+1})
{\rm \ is\ nonvanishing} \ \forall\ 1\leq i\leq n
\end{equation}
where we use the convention that for indices $n+1=1$.

\smallskip

{\it How can we determine the existence of such surfaces?}

\smallskip

Were we looking for $n$-folds, the answer
would be given by the Frobenius theorem:

\begin{theorem}[Frobenius theorem] 
Given pointwise linearly independent one-forms $\psi^1\hd \psi^n$ on a
manifold $X^m$, there exists an immersed
submanifold $i:M^{m-n}\ra X$ passing through
$p\in X$ on which $i^*(\psi^j)\equiv 0$ for all $j$ 
(i.e.,  with $T_xM=\Delta_x:=\tker\{\psi_j\}$ for all $x\in M$)
if in a neighborhood of $p$ there
exist one-forms $\a^i_j\in \Omega^1(X)$ such that
\begin{equation}\label{frob1}
d\psi^i=\a^i_1\ww\psi^1+\cdots+ \a^i_n\ww\psi^n \ \forall \ i
\end{equation}
\end{theorem}

The condition \eqref{frob1} is often expressed
as $d\psi^i\equiv 0\tmod \{ \psi^1\hd \psi^n\}$.
In fact the individual forms don't matter, just
their span, so we could   write
$I=\{ \psi^1\hd \psi^n\}$ and 
$$d\psi^i\equiv 0\tmod I \ \forall \ i.
$$
Another way to express it is that
locally if, $X,Y$ are vector fields lying in $\Delta$,
that $[X,Y]$ also lies in $\Delta$. (Exercise: verify
that this is indeed equivalent.)  Note that
all these conditions involve beginning with
first order information and differentiating it
once - if everything is OK, then we are guaranteed
solutions. That is, we can stop working after
taking two derivatives.

Were we in the situation that there was just a single
one-form, then {\it Pfaff's theorem} (see, e.g., \cite{BCG3} \S 1.3)
guarantees 
existence of submanifolds of dimension roughly half the
dimension of the manifold. Moreover, by computing the
exterior derivative of the one-form one can
determine the precise maximal dimension of
a submanifold on which the form pulls back to be zero.

To deal with the general setting of determining existence
of submanifolds on which an ideal of differential forms pulls
back to be zero, an explicit algorithm
was developed by Cartan and others. The algorithm also gives
a rough estimate of the size of the space of such
manifolds. (E.g., in the Frobenius theorem, there
is a unique such manifold through a point but for
Pfaff's theorem, there will be \lq\lq functions \rq\rq\
worth of solution manifolds through a point.)

The essential question is: Given a candidate tangent
space (a first order admissible Taylor series), can
we extend it? - i.e., can we \lq\lq fit together\rq\rq\ 
potential tangent spaces to obtain a solution submanifold?

\section{EDS terminology}
Let $V$ be a vector space,
let $G(k,V)$ denote the {\it Grassmannian} of $k$-planes
through the origin in $V$.

\begin{definition}
Let $\Sigma$ be a manifold
Let $\ci\subset \Omega^*(\Sigma)$ be a differential
ideal, which we will call an {\it exterior differential system}.
We let $\ci_j\subset \Omega^j(\Sigma)$ denote the component in
degree $j$ and
we will henceforth assume $\ci_0=\emptyset$.
An {\it integral manifold} of $\ci$ is an immersed
submanifold $i:M\ra \Sigma$ such that $i^*(\phi)=0$
for all $\phi\in\ci$.
\end{definition}

 As with many things in
mathematics, we will work infinitesimally with the
goal  of linearizing the problem of determing
the integral manifolds of an EDS.

\begin{definition}
For $x\in \Sigma$, we let  
$$\cV_k(\ci)_{x}:=\{ E\in  G(k, T_x\Sigma) \mid
\phi|_E=0\ \forall\phi\in\ci\}
$$
which is called the variety of {\it $k$-dimensional integral
elements} to $\ci$ at $x$. We let
$\bold G(k,T\Sigma)$ denote the Grassmann bundle, i.e.,
the bundle over $\Sigma$ whose fiber over $x\in \Sigma$
is $G(k,T_x\Sigma)$,
and let $\cV_k(\ci)\subset \bold G(k,T\Sigma)$ denote
the set of all $k$-dimensional integral elements.
\end{definition}

The first step in the Cartan algorithm is {\it Cartan's test}:
one compares a crude estimate (obtained from linear algebra
calculations) of $\tdim\cV_k(\ci)$ with its actual
dimension. If the two numbers agree, then the
{\it Cartan-K\"ahler theorem} guarantees local existence
of integral manifolds.
We can think of it as saying \lq\lq if the second order
terms for the Taylor series look good, everything is good\rq\rq .
If the test fails, we must take more derivatives to 
determine existence. This process is called {\it prolongation}.
The {\it Kuranishi prolongation theorem} says
that in principle one only needs to prolong a finite
number of times before getting a definitive answer, but
this is of little use in practice as the theorem
gives no indication of how many times one must prolong
(how many derivatives one needs to take).
Before going into details, let's examine some
examples to develop our intuition.

\section{PDE and EDS}
\begin{example}\label{otnine}
Consider the PDE system for $u(x,y)$ given by
\begin{equation}\label{basicsys}
\begin{aligned}
u_x &= A(x,y,u), \\ u_y &= B(x,y,u),
\end{aligned}
\end{equation}
where $A,B$ are given smooth functions.
Since \eqref{basicsys} specifies both partial derivatives of $u$,
at any given point $p=(x,y,u)\in \BR^3$
the tangent plane to the graph of a solution passing through $p$
is uniquely determined.

Whether or not the plane is actually tangent to
a solution to \eqref{basicsys} depends
on whether or not the equations \eqref{basicsys}
are \lq\lq compatible\rq\rq\
as differential equations.
For smooth solutions to a system of PDE,
compatibility conditions arise because mixed partials
must  commute, i.e.,   $(u_x)_y = (u_y)_x$.  In our example,
$$
\begin{aligned}
(u_x)_y  & = \frp  {y}   A(x,y,u) = A_y(x,y,u) + A_u(x,y,u)\frac{\partial u}{\partial y}  = A_y +BA_u,
\\
(u_y)_x & = B_x +AB_u, \end{aligned}
$$
so setting $(u_x)_y = (u_y)_x$ reveals a ``hidden equation'', the compatibility condition
\begin{equation}\label{hid}
A_y + BA_u = B_x + AB_u.
\end{equation}
and the Frobenius condition is exactly the vanishing
of this equation. To see this let
$$
\th= du- A(x,y,u)dx -B(x,y,u)dy.
$$
Exercise: show that \eqref{hid} holds iff 
$d\th\equiv 0\tmod \th$.

Here we have the EDS $\ci=\{\th\}_{diff}$ on $\Sigma=\BR^3$ but since this EDS comes
from a PDE, we have an additional condition that we want our integral
manifolds to satisfy, namely that $x,y$ are independent variables on
a solution. We encode this by setting $\Omega= dx\ww dy$ and making
the following definitions:

\begin{definition}
Let $\ci\subset \Omega^*(\Sigma)$ be a differential
ideal, and $\Omega\in \Omega^n(\Sigma)$. The pair
$(\ci,\Omega)$ is called   an {\it exterior differential system
with independence condition}.
An {\it integral manifold} of $\ci$ is an immersed
submanifold $i:M\ra \Sigma$ such that $i^*(\phi)=0$
for all $\phi\in\ci$ and $i^*(\Omega)$ is nonvanishing.
Note that we really only need $\Omega$ up to scale and modulo
$\ci$, so we sometimes refer to an independence condition
as an equivalence class of $n$-forms (the equivalence is up to scale and
modulo $\ci$).
\end{definition}

\begin{remark} One can attempt to obtain
solutions to the system \eqref{basicsys}
by solving a succession of Cauchy problems. For
example fix $y= 0$ and solve the ODE
\begin{equation}\label{odea}
\dfrac{d\tilde u}{dx} = A(x,0,\tilde u),\qquad  \tilde u(0) =u_0.
\end{equation}
After solving \eqref{odea}, hold $x$ fixed and solve   the initial
value problem
\begin{equation}\label{odeb}
 \dfrac{du}{dy} = B(x,y,u),\qquad  u(x,0) = \tilde u(x).
\end{equation}
This determines a function $u(x,y)$ on some neighborhood of $(0,0)$.
The problem is that this function may not satisfy our original
equation, and it also may depend
on the path chosen. The function is independent of path chosen precisely
if the Frobenius condition holds, and in that
case it gives the right answer too.
\end{remark}
\end{example}

\def\FF#1#2{f^{#1}_{#2}}
\def\uu#1#2{u^{#1}_{#2}}

In general,
given a first-order system of $r$ equations  for $s$ functions $u^a$ of
$n$ variables, there exists a change of coordinates so that
the system takes the form
$$
\begin{aligned}
&\uu 1{\xx 1}= \FF 11({\overline x},{\overline u}),\\
&\vdots \\
&\uu{r_1}{\xx 1} = \FF1{r_1} ({\overline x},{\overline u}),\\
&\uu 1{\xx 2} = \FF 21  ({\overline x},{\overline u}, {\overline u}_{\xx 1}),\\
&\vdots \\
&\uu {r_2}{\xx 2} = \FF 2{r_2} ({\overline x},{\overline u}, {\overline u}_{\xx 1}),\\
&\vdots\\
&\uu 1{\xx n} = \FF n1  ({\overline x},{\overline u}, {\overline u}_{\xx 1}\hd{\overline u}_{\xx{n-1}}),\\
&\vdots \\
&\uu {r_n}{\xx n} = \FF n{r_n} ({\overline x},{\overline u}, {\overline u}_{\xx 1} \hd  {\overline u}_{\xx{n-1}}),
\end{aligned}
$$
where $\overline x=(x^1\hd x^n)$, $\overline u= (u^1\hd u^s)$,
$u^a_{\xx j} = \frpp{u^a}{\xx j}$,
$1\leq a\leq s$, $1\leq  j\leq n$, and $r_1\leq r_2\leq\hdots\leq r_n=s$ with $r= r_1+\hdots + r_n$
(see \cite{cartan1931}).

We may be able to produce solutions of this system by solving a
series of Cauchy problems.  However, we need to check that
equations are compatible, i.e.,
that mixed partials commute:
$$
\frpp{}{\xx  i}\FF  j{\s} = \frpp{}{\xx j}\FF i{\s},\qquad
1\leq i,j\leq n,\ \forall \s.
$$
Although it would be impractical to change any given system of PDE into the above form, converting this system to an EDS will guide  us naturally
to the analog of the above form.  We can then apply a straightforward test
that signals when no further compatibility conditions need to be checked. 

\section{Cartan's test}
Let $\ci$ be an EDS on a manifold $\Sigma$. 
Let $p\in \Sigma$ be a general point and
$(p,E)\in \cV_{n}(\ci)$ be a general point
of $\cV_{n}(\ci)$.  The required generality can
be made precise, see \cite{IvL,BCG3} but we suppress
that in these lectures. Intuitively, we want 
$(p,E)$ to be \lq\lq like\rq\rq\ its neighbors in
some small open set in $\cV_n(\ci)$.

\begin{remark}
Note that since we are dealing
with (analytic) varieties, i.e.,  zero sets of  
analytic functions, there can be components
to $\cV_{n}(\ci)_{p}$. \lq \lq A general point\rq\rq\  means a general point 
of a given component.
\end{remark}

As mentioned above, the test we are after will compare
a codimension estimate obtained by linear algebra
calculations with the codimension of
a variety.  

\begin{definition} Let $E\in \cV_{j}(\ci)_{p}$
and let $e_1\hd e_j$ be a basis of $E$. Define
$$
H(E):=\{ v\in T_p\Sigma\mid \phi(v,e_1\hd e_j)=0\ 
\forall \phi\in \ci_{j+1}\}
$$
the {\it polar space} of $E$.
\end{definition}

 Note that
\begin{enumerate}
\item   $H(E)$ is well defined (i.e., independent
of our choice of basis), 
\item  $E\subset H(E)$ and 
\item  determining
$H(E)$ is a linear calculation.
\end{enumerate}

The quotient $H(E)/E$ may be thought of as the space of
possible enlargements of $E$ from a $p$-dimensional
integral element to a $(p+1)$-dimensional
integral element.  We will actually need to
calculate the dimensions of a series of polar spaces.

Let $E\in \cV_n(\ci)$. Fix a generic
flag $E_1\subset\cdots \subset E_{n-1}\subset E_n=E$
in $E$. Let $c_j=\tcodim(H(E_j), T_p\Sigma)$,
and set $c_0=\tcodim \cV_1(\ci)=\tdim \ci_1$. Note
that if $\Sigma$ has components, then $\tcodim \cV_1(\ci)$
can depend on the component, and for $j>1$, $\cV_j(\ci)$
may have components even if $\Sigma$ has just one component.
Therefore we will write $\tcodim_{E_{j+1}}(\cV_{j+1}(\ci),\bold G(j+1,T\Sigma))$
to eliminate any possible ambiguity when discussing
the codimension of $\cV_{j+1}(\ci)$ at $E_{j+1}$.

We have the following estimate
\begin{proposition} 
$$
\tcodim_{E_{j+1}}(\cV_{j+1},\bold G(j+1,T\Sigma))
\geq \tcodim_{E_j}(\cV_j,\bold G(j,T\Sigma)) +
\tcodim_{T_p\Sigma}H(E_j).
$$
\end{proposition}

The inequality is intuitively reasonable as the first term on the right
represent the conditions to have a $j$-dimensional integral 
element and the second term represents the new conditions for
enlarging it to a $(j+1)$-dimensional integral element.  
Equality holding should be interpreted as
$\cV_{j+1}$ being \lq\lq as large as possible\rq\rq\  
at $E_{j+1}$. Adding up these inequalities, we obtain

\begin{proposition}
\begin{equation}\label{cartinequal}
\tcodim_E(\cV_n, \bold G(n, T\Sigma)\geq c_0+c_1+\cdots +c_{n-1}.
\end{equation}
\end{proposition}

The {\it Cartan-K\"ahler theorem} states that when equality holds
(assuming our genericity hypotheses about $p$ and $E$),
there exists an $n$-dimensional integral manifold
through $p$ with tangent space $M$.
The test  for equality holding in \eqref{cartinequal} is called
{\it Cartan's test}. If an integral element passes Cartan's
test, we get a bonus - a coarse estimate of the size  of
the moduli space of integral manifolds through $p$. Namely if we set $s_k=c_k-c_{k-1}$ and
let $k_0$ be the largest integer
such that $s_{k_0}$ is nonzero, then integral
manifolds depend roughly on $s_{k_0}$ analytic functions
of $k_0$ variables. In particular if the largest $k_0$ is
$0$, then integral manifolds depend only on a choice of
constants, as in the Frobenius theorem.

\subsection*{Other possibilities.}

\subsubsection*{$\cV_{n }(\ci)_p=\emptyset$.}
More precisely, there exists a Zariski open subset
of $\Sigma$ over which there are no $n$-dimensional
integral elements. In this case it is necessary to
restrict to the (analytic) subvariety $\Sigma'\subset \Sigma$
over which there are $n$-dimensional integral elements and start over,
working at general points of $\Sigma '$. Note that
$\Sigma '$ may have several components and that one
must perform the test on each component separately.
If $\tdim \Sigma' <n$ we are done, there are no
$n$-dimensional integral manifolds.

\subsubsection*{{\it Cartan's test fails.}}
Intuitively, this means we have not differentiated enough
to uncover all compatibility conditions and we must
take more derivatives. It turns out 
that, rather than taking higher
derivatives, it is notationally simpler   to start over on a larger space where
our old derivatives are replaced by independent variables.
(This corresponds to the standard process of converting
any system of PDE to a first order system by adding 
additional variables.)

More precisely, forgetting about $\ci$ for the moment, on
$\pi: \bold G(n,T\Sigma) \ra \Sigma$, consider the following
tautological system:
given $(p,E) \in \bold G(n,T\Sigma)$, we
have $E\upperp \subset T^*_p\Sigma$. Define
$$
I_{(p,E)}:=\pi^*(E\upperp).
$$
For good measure we add
the independence condition determined by
$\Lambda^n(\pi^*(T^*\Sigma/I))$. Integral
manifolds of the tautological system $(\{ I\}_{diff}, \Omega)$
with $[\Omega]\in \Lambda^n(\pi^*(T^*\Sigma/I))$
are precisely the {\it Gauss images}
of immersed $n$-dimensional submanifolds
$f: M\ra \Sigma$.

Now let's return to our original EDS $\ci$ on $\Sigma$:

\begin{definition} The {\it prolongation} of $\ci$ is
the pullback of the tautological system on $\bold G(n,T\Sigma)$
to $\cV_n(\ci)\subset \bold G(n,T\Sigma)$.
\end{definition}

One then starts over with $\Sigma$ replaced
by $\cV_n(\ci)$ and $\ci$ replaced by
  the pullback of the tautological system. One
then performs Cartan's test, if it fails, one
prolongs again, etc... For more details, see
\cite{IvL}, \S 5.5.

\section{First examples of Cartan's test}

\subsection*{Example 0: arbitrary maps $\BR^2\ra \BR^2$.}
Let the first $\BR^2$ have coordinates $x^1,x^2$,
the second coordinates $u^1,u^2$ and
let $\Sigma= J^1(\BR^2,\BR^2)\simeq \BR^8$ with
coordinates $(x^i,u^j,p^i_j)$, $1\leq i,j\leq 2$.

Given a map $f:\BR^2\ra\BR^2$, we define the
{\it lift} of $f$ to $\Sigma$ to be the
set of points 
$$(x^1,x^2,f^1(x),f^2(x),\frpp {f^1}{x^1}|_x
\frpp {f^1}{x^2}|_x,\frpp {f^2}{x^1}|_x,
\frpp {f^1}{x^2}|_x),
$$ 
which is a coordinate version of the Gauss map of
an immersion.
Let
\begin{align*}\ta 1&= du^1-p^1_1dx^1-p^1_2dx^2\\
\ta 2&= du^2-p^2_1dx^1-p^2_2dx^2
\end{align*}
Introduce the independence condition $\Omega=dx^1\ww dx^2$.
Then integral manifolds of the system
$(\{\ta 1,\ta 2\}_{diff},\Omega)$ are in one to one
correspondence with lifts of maps $f:\BR^2\ra\BR^2$.

The manifold $\Sigma=J^1(\BR^2,\BR^2)$ equipped with the
system $(\{\ta 1,\ta 2\}_{diff},\Omega)$ is  
called the space of {\it one-jets} of
mappings $\BR^2\ra \BR^2$.

Let's perform Cartan's test:

\subsection*{Determination of $c_0+c_1$.}
$c_0=2$ because $\tdim\ci_1=2$. The equations
on any line $\{ v\}$ are explicitly
$\ta 1(v)=0$, $\ta 2(v)=0$.

To find $c_1$, we need to take a generic $\{ v\}\in \cV_1$.
Write
$$
v= a^1\frp { x^1}+a^2\frp { x^2}+b^i_a \frp {p^a_i}
+e_a \frp {u^a}
$$
where here and throughout we use the {\it summation convention}
that repeated indices are to be summed over.
$\ta j(v)=e_j-p^j_1a^1-p^j_2a^2$ so we may 
take
$$
v= a^1\frp { x^1}+a^2\frp {x^1}+b^i_a \frp {p^a_i}
+(p^a_1a^1+p^a_2a^2)\frp {u^a}
$$
where $a^j,b^i_a$ are (general) constants.

To determine $c_2$, we must find $\ci_2$.
First there are $\a\ww \ta 1,\a\ww \ta 2$ where
$\a$ is any one-form. We also have
$d\ta j= -dp^j_1\ww dx^1-dp^j_2\ww dx^2$.
To determine a possible enlargement of
$\{ v\}$ we must calculate
$$
d\ta j(v,\cdot)=b^j_1 dx^1 -a^1dp^j_1+b^j_2dx^2-a^2dp^j_2
$$ 
So any vector $w$ in $H^1(\{v\})$ must
satisfy the four linear equations
$$
\ta j(w)=0, d\ta j(v,w)=0
$$ 
These are independent (check yourself!),
so we
obtain $c_1=4$ and $c_0+c_1=6$.

\subsection*{Determination of $\tcodim \cV_{2}$.}
Let $\bold G(2,T\Sigma)$ have
local coordinates $(x^i,u^a,p^a_i; b^i_a,c^i_a,e^a,f^a)$
where the first set gives coordinates for the base point
and the second for the plane $v\ww w$ where
\begin{align*}
v&= \frp{x^1}+b^i_a\frp{p^a_i}+ e^a\frp{u^a}\\
w&= \frp{x^2}+c^i_a\frp{p^a_i}+ f^a\frp{u^a}
\end{align*}
We have the following conditions and consequences:
$$
\begin{array}{ccc}
\ta 1(v)=0 &\implies & b^1=p^1_1\\
\ta 2(v)=0 &\implies & b^2=p^2_1\\
\ta 1(w)=0 &\implies & c^1=p^1_2\\
\ta 2(w)=0 &\implies & c^2=p^2_2\\
d\ta 1(v\ww w)=0&\implies & c^1_1-b^1_2=0\\
d\ta 2(v\ww w)=0&\implies & c^2_1-b^2_2=0
\end{array}
$$
These six equations are independent and
we conclude $\tcodim \cV_2(\ci)=6$ and
Cartan's test succeeds. Moreover integral manifolds
\lq\lq depend on  two functions of two variables\rq\rq\
which in this case we see explicitly, as we knew the
solutions all along.

\subsection*{Example 1: The Cauchy-Riemann equations
$u^1_{x^1}=u^2_{x_2},u^1_{x^2}=-u^2_{x_1}$.}
This example is the same as above, except that
we now restrict to the submanifold $\Sigma'\subset \Sigma$
where $p^1_1=p^2_2$ and $p^1_2=-p^2_1$.
We still have $c_0=2$ as $\ta 1,\ta 2$ remain
linearly independent when restricted to $\Sigma'$
but we now have $c_1=2$ (exercise - be sure to
express the initial $v$ in terms of $6$ variables
(e.g., eliminate $p^2_1,p^2_2$)).
Similarly, only four of the six equations for
$\cV_2$ remain independent.
So here we have the equality $\tcodim_E\cV_2=
4=c_0+c_1=2+2$.

Here Cartan's test indicates that integral manifolds
should depend on two functions of one variable, which
we also know to be the case as a (sufficiently generic)
real analytic arc uniquely determines a holomorphic map $\BC\ra \BC$.

\begin{remark} Note that in both the above calculations,
the calculation of $\tcodim \cV_2(\ci)_p$ was linear.
There is a large class of EDS, called {\it linear Pfaffian systems}
which are systems defined by one-forms for which this linearity
holds. For such systems, there
is a simplified version of Cartan's test. Any system of partial differential
equations expressed as the pullback of the contact
system on the space of jets is a linear Pfaffian system,
see, e.g., \cite{IvL}, example 5.1.4.
\end{remark}

\subsection*{Example 2: Lagrangian   submanifolds.}
\label{Lagex}
Let $\omega$ be the standard symplectic form on $\BR^{2n}$:
\index{symplectic form}
$$\omega= dx^1 \wedge dy^1 + \hdots + dx^n \wedge dy^n.$$
An $n$-dimensional submanifold is {\em Lagrangian}
\index{submanifold!Lagrangian} if it is an integral manifold
of $\ci=\{\omega\}_{diff}$.

Given $(p,E)\in \cV_n(\ci)$, we can make a linear change of coordinates (while keeping
the form of $\omega$) so that $E$ is annihilated by $dy^1\hd dy^n$.
This is because 
the subgroup of $GL(T_p\BR^{2n})$ leaving $\omega$ invariant is the symplectic group which
acts transitively not only on Lagrangian $n$-planes but
on all flags within them. 
Thus all $n$-planes at all points are equivalent and genericity
issues don't enter.
 Any nearby integral
$n$-planes at $p$ are given by $dy^j=\sum_k s^{jk} dx^j$ for $s^{jk}=s^{kj}$.  Therefore,   $\tdim(\cV_n(\ci)_p)=\binom{n+1}2$, 
$$\tcodim_E(\cV_n(\ci)_p, G(n,T_p\BR^{2n}))=
\tcodim_E(\cV_n(\ci), {\bold G}(n,T\BR^{2n}))=\binom n2,
$$
 independent
of $p$ and $E$.

Let $e_1\hd e_n \in E$ be dual to $dx^1\hd dx^n$
and we use $e_1\hd e_n$ to build our flag in $E$, i.e.,
$E_j=\langle e_1\hd e_j\rangle$.  
(By the remark above, there are no genericity issues to be concerned with.)

It is easy to calculate that for $j\leq n$,
$$
H(E_j)=\{ v\in T_p\BR^{2n}\mid dy^k(v)=0\ \forall \ k\leq j\}
$$
so $c_j=j$ for $j\leq n-1$.  Since $c_1 + c_2 + \hdots +c_{n-1} = \binom{n}2$,
we have involutivity, and integral manifolds depend on 1 function of $n$ variables.
(In fact, they can be explicitly constructed by setting $y^j=\partial f/\partial x^j$ for
$f$ an arbitrary function of $x^1\hd x^n$.)

\section{Periodic billiard orbits}

We now return to the problem of finding $n$-periodic billiard orbits.
(The results of this section and the next are joint with Baryshnikov and Zharnitsky.)
We have the EDS
$$
\ci=\{ \psi_i\}_{diff}
$$
and several independence conditions: that each $\eta^i\ww \eta^{i+1}$
is nonvanishing on an integral manifold. Fortunately we can
reduce to a single independence condition thanks to the following
lemma:

\begin{lemma} It is sufficient to work with the independence
condition $\eta^1\ww \eta^2$ (or any $\eta^i\ww\eta^{i+1}$).
\end{lemma}

\begin{proof} 
Let $X_i$ be a dual basis to $\eta^j$ of $\tker\ci_1$.
Take local coordinates $p^{\a}_1,p^{\a}_2$ about
$[X_1\ww X_2]$ where we write $[v\ww w]$ as a nearby point with
$v= X_1 + \ppp\a 1X_{\a}$, $w=X_2+\ppp\a 2X_{\a}$.

Introduce the notations $a_j=\frac{\tcos (\a_{j+1})}{2l_j}$,
$b_j=\frac{\tcos (\a_{j-1})}{2l_{j-1}}$ where we use the notation
of \S\ref{setsect}. One calculates
(see \cite{BLZ}) that 
$$
d\psi^j\equiv (a_j\eta^{j+1}+b_j\eta^{j-1})\ww\eta^j \tmod I
$$
Moreover $p\in \Sigma$ implies that none of the $a_j,b_j$ are zero
at $p$.

Evaluating the $d\psi^i$ at $v\ww w$ (that is, evaluating
at an arbitrary point in our chart) we obtain the $n$  equations
on the $p^{\a}_j$
\begin{equation}\label{pajeqns}
\begin{array}{ccc}
0&=& a_1 + b_1\ppp n2\\
0&=&\ppp 31a_2+b_2\\
0&=& (\ppp 42\ppp 31 -\ppp 32\ppp 41)a_3 +\ppp 31b_3\\
0&=&  (\ppp 52\ppp 41-\ppp 51\ppp 42)a_4+(\ppp 42\ppp 31 -\ppp 32\ppp 41)b_4\\
&&\vdots \\
0&=& \ppp n2a_n +  (\ppp n2\ppp{n-1}1-\ppp{n-1}2\ppp n1)b_n
\end{array}
\end{equation}
of which $n-1$ are independent.

The first equation implies $\ppp n2\neq 0$, which implies that
on an integral element on which $\eta^1\ww \eta^2\neq 0$,
we also have $\eta^n\ww\eta^1\neq 0$. The second
equation implies $\ppp 31\neq 0$ which implies that  similarly
$\eta^2\ww\eta^3\neq 0$. The third equation implies
that $(\ppp 42\ppp 31 -\ppp 32\ppp 41)\neq 0$ but this
is exactly the condition that $\eta^3\ww\eta^4\neq 0$. Continuing,
we see that $\eta^1\ww\eta^2\neq 0$ implies that
all $\eta^j\ww\eta^{j+1}\neq 0$ on an integral element.
\end{proof}

\begin{remark} Had we instead taken, e.g., $\eta^1\ww\eta^3$ as independence
condition, (assuming $n>3$) we could not have drawn a similar conclusion,
see \cite{BLZ}.
\end{remark}

Introduce notation $\Delta_{j-1}=(\ppp j1\ppp{j-1}2-\ppp j2\ppp{j-1}1)$ with the
convention that $\ppp 11=\ppp 22=1$, $\ppp 21=\ppp 12=0$, so $\Delta_1=1$.
Then our equations \eqref{pajeqns} become
$$
a_j\Delta_{j-1}+b_j\Delta_j=0
$$
which we may write in matrix form:
$$
\begin{pmatrix}
0&0&0&\hdots&0&a_1\\
a_2&b_2&0&\hdots&0&0\\
0&a_3&b_3&\hdots&0&0\\
\vdots&\vdots&\vdots&\vdots&\vdots&\vdots\\
0&0&0&\hdots&a_{n-1}&b_{n-1}\end{pmatrix}
\begin{pmatrix}
\Delta_2\\ \Delta_3\\ \Delta_4\\ \vdots \\ \Delta_n\end{pmatrix}
=
\begin{pmatrix} -b_1\\ -a_2\\ 0\\ \vdots \\ 0\end{pmatrix}
$$
and since the $a_j$ are nonzero, there is a unique solution for
$\Delta_2\hd \Delta_n$.
Now $\Delta_2=\ppp 31$ and $\Delta_n=\ppp n2$ so $p^3_1,p^n_2$ are   fixed
and the remaining equations on the $\ppp\a i$
  are
independent. In fact one can solve explicitly for all the
remaining $\ppp \a 1,\ppp\a 2$  in terms of $\ppp 41,\ppp 32,\ppp 51,\ppp 61\hd \ppp{n-1}1$.
Thus the space of integral elements satisfying the genericity condition
is of dimension $n-3$.

\begin{proposition} The system $(\ci,\eta^1\ww \eta^2)$ has $\tcodim \cV_2(\ci )= 3n-1$,
$c_0=n$, $c_1=2n-2$ and thus fails Cartan's test by one.
\end{proposition}

\begin{proof}
Here $c_0$ is just the codimension of the space of one-dimensional
integral elements at a point of $\Sigma$. To calculate $c_1$,
one needs a sufficiently generic vector, $Z=X_1+\cdots + X_n$ will
do. One then sees that $Z$ is contained in a unique two-dimensional
integral element.
\end{proof}

If one ignores the genericity conditions,
as $n$ increases the dimensions
of integral manifolds can be   arbitrarily large  
(see \cite{BLZ}).
The next proposition states that with the genericity conditions,
this fails even at the infinitesimal level.

\begin{proposition} For all $n$, there are no $3$-dimensional integral
elements to $\cI$ satisfying the genericity conditions.
\end{proposition}
\begin{proof}
On a three dimensional integral element, we must have
say $\eta^1, \eta^2, f_a\eta^a$ independent where $3\leq a\leq n$
and the $f_a$'s are some constants.
First note that $f_3,f_n$ must be zero by considering
$d\psi^2$, $d\psi^1$ respectively. But we also must have
$\eta^2, \eta^3$ independent, and since $f_3=0$, this implies
$\eta^2,\eta^3, f_a\eta^a$ must be independent, which, using
$d\psi^3$ implies that $f_4=0$. Continuing in this fashion
one obtains that all the $f_a$ must be zero.
\end{proof}

\section{Three periodic billiard orbits}
Here the space of integral elements satisfying the billiard conditions is
a single point. Taking $\eta^1\ww\eta^3$ as our independence condition,
writing $c_j=\tcos(\a_j)$, $s_j=\tsin(\a_j)$,
we see that on integral elements
\begin{equation}\label{n3eqn}
\eta^2+\frac{c_1l_2}{c_2l_3}\eta^1+\frac{c_3l_1}{c_2l_3}\eta^3=0
\end{equation}
Adding  this form to the ideal and taking
its derivative, we see
\begin{align*}
 &d(\eta^2+\frac{c_1l_2}{c_2l_3}\eta^1+\frac{c_3l_1}{c_2l_3}\eta^3)\\ 
 &\  \equiv  
[(-s_3c_1c_2+c_3s_2c_1+c_3s_1c_2)l_1
+(-c_3s_1c_2+s_3c_1c_2+c_3s_2c_1)l_2
\\
&\ \   
+(-c_3s_2c_1+s_3c_1c_2+c_3s_1c_2)l_3]\frac{\eta^1\ww\eta^3}{c_2^2l_3^2}
\end{align*}
Thus $\cV_1(\ci)_x=\emptyset$ for general $x\in \Sigma$ and
we must restrict to the subvariety of $\Sigma$ where
\begin{align*}
&(-s_3c_1c_2+c_3s_2c_1+c_3s_1c_2)l_1
+(-c_3s_1c_2+s_3c_1c_2+c_3s_2c_1)l_2\\
&
+(-c_3s_2c_1+s_3c_1c_2+c_3s_1c_2)l_3=0.
\end{align*}
Now recall that a triangle is uniquely determined, e.g., by two
of its three angles and the length of one of its sides, we may
write
\begin{align*}
\a_3&=\frac{\pi}2 -\a_1-\a_2\\
l_3&=\frac{l_1\tsin (2\a_2)}{\tsin(\pi -2\a_1 -2\a_2)}\\
l_2&=\frac{l_1\tsin (2\a_1)}{\tsin(\pi -2\a_1 -2\a_2)}\end{align*}
and substituting in, we obtain the equation
$$
6l_1c_1c_2s_1s_2=0
$$
which cannot occur on $\Sigma$.\qed

\section{A few successes of the Cartan-K\"ahler theorem}

\subsection{The Cartan-Janet theorem}
Given an analytic Riemannian manifold ($M^n,g)$, does
there exist a local isometric immersion into
Euclidean space $\BE^{n+s}$? The Cartan-Janet theorem
states that for any analytic metric the answer
is yes as long as $s\geq \binom n2$. If the metric
is special one can sometimes do much better, see
\cite{C1919,BBG,BGY} for the cases of
space forms and generalizations.

\subsection{Manifolds with exceptional holonomy}
Using EDS 
Bryant \cite{bryantg2} showed that there exist
non-symmetric
Riemannian manifolds with holonomy $G_2$ and
$Spin_7$, settling the last open local existence
questions in the Riemannian case of Berger's 
1953 thesis \cite{berger}.

\subsection{Existence of calibrated submanifolds}
The abundance of special Lagrangian and other calibrated
submanifolds was first proved by Harvey and Lawson
\cite{HL} using the Cartan-K\"ahler theorem.

While describing the first two examples would involve too many
definitions, we will explicitly describe two cases of
applying the Cartan-K\"ahler theorem to prove
existence of calibrated submanifolds.

\begin{definition} A {\it calibration} on an
oriented Riemannian manifold $\Sigma$ is
a closed differential form $\phi\in \Omega^k(\Sigma)$ such that
for all unit volume $(p,E)\in \bold G(k, T\Sigma)$,
$\phi(E)\leq 1$.
\end{definition}

 There are many variants on the definition. Calibrations
are a tool for finding volume minimizing submanifolds of $\Sigma$ 
because the fundamental lemma of calibrations says that
if $i:M\ra \Sigma$ is an immersed submanifold on which
$i^*(\phi)=\tvol_M$ then $M$ is volume minimizing in its
homology class (assuming $M$ is compact, there are
variations when $M$ is noncompact), see \cite{HL}.

Recently calibrated manifolds have become of central importance
because of applications to physics. See, e.g., 
Joyce's lectures in \cite{joyce}. 
Calibrations may be thought of as generalizations
of normalized powers of the K\"ahler form, which itself
gives rise to an involutive system (the Cauchy-Riemann equations!).
We will discuss two additional calibrations, the {\it special
Lagrangian calibration} and the {\it associative calibration}.

Sometimes a calibration
$\alpha$   has a complementary form $\alpha_c$ such that
\begin{equation}\label{comp}
|\alpha(E)|^2+ |\alpha_c(E)|^2=1
\end{equation}
for all unit volume planes $E$.
In such cases we may define an EDS whose integral
manifolds are the submanifolds calibrated by $\a$ by 
taking $\ci=\{ \a_c\}_{diff}$.
 
\begin{example}[Special Lagrangian manifolds]
On $\BR^{2n}=\BC^n$ (or
any K\"ahler manifold), consider the differential $n$-form
$$\alpha = \realpart\left(dz^1 \wedge \cdots \wedge dz^n\right),$$
where $z^j=dx^j+i dy^j$, called the {\it special Lagrangian calibration}.

  In the special Lagrangian case, a variant of \eqref{comp} holds. If
we take
$$\a_c = \impart\left(dz^1 \wedge \cdots \wedge dz^n\right),$$
then, restricted to {\it Lagrangian $n$-planes}, \eqref{comp}
holds.
Moreover, it is easy to see any submanifold calibrated
by $\a$ is Lagrangian, so   $\ci=\{ \omega,\a_c\}_{diff}$
is an EDS whose integral manifolds are the special Lagrangian
submanifolds.

Given $E\in \cV_n(\ci)$, we can   change coordinates so that
$E$ is annihilated by $dy^1\hd dy^n$. (This is because
the system is $SU(n)$ invariant and $SU(n)$ acts transitively
on the special Lagrangian planes at a point and even transitively
on flags in special Lagrangian planes.)  Taking
$e_1\hd e_n \in E$ to be dual to $dx^1\hd dx^n$, we have
$c_j=j$ for $1\leq j\leq n-2$ as in example 2 of \S\ref{Lagex}.  However,
$$\left.
\begin{aligned} \omega &\equiv dx^{n-1} \wedge dy^{n-1} + dx^n \wedge dy^n \\
\alpha_c &\equiv dx^1 \wedge \cdots dx^{n-2} \wedge (dx^{n-1} \wedge dy^n - dx^n \wedge dy^{n-1})
\end{aligned}\right\} \tmod  dy^1\hd dy^{n-2}
$$
shows that $c_{n-1}=n$.  The requirement that $\alpha_c|_{E}=0$
is one additional equation
($\sum_j s^{jj}=0$) on the set of Lagrangian $n$-planes 
so   the codimension of $\cV_{n}(\ci)$ is one greater than the Lagrangian case and the system is   involutive, with solutions
depending on two functions of $n-1$ variables.
\end{example}

\begin{example}[Associative submanifolds]
\label{seveng2ex}
\index{submanifold!associative}
\index{associative submanifolds}
The 14-dimensional compact Lie group $G_2$ arises as the automorphism group of
the normed algebra $\BO$ of octonions (see e.g., \cite{IvL} \S A.5)), and leaves invariant
a 3-form $\phi$ on $\BR^7=\impart \BO$,
where $\phi(x,y,z)=\langle x,yz\rangle$. (Here $\langle\cdot,\cdot\rangle$
is the inner-product induced from the norm.)   This $\phi$  is
a calibration on $\BR^7$, and it admits a complement as in \eqref{comp}:
$\phi_c=\tfrac12  \impart\left((xy)z-(zy)x\right) $.

We define an EDS $\ci$ for associative submanifolds by taking
the components of the $\impart\BO$-valued 3-form $\phi_c$ as generators.
(Since $\phi_c$ is constant-coefficient, all of these generators are closed.)

Let $E\in \cV_3(\ci)$. Then the stabilizer of $E$ in
$G_2$ is six-dimensional.  Since   $G_2$ acts transitively on the space
of $3$-dimensional integral elements, we conclude
$$ \tcodim (\cV_3(\ci)_p,G(3,T_p\impart \BO))=12-8=4.$$  
On the other hand, for any
flag in $E$, $c_0=c_1=0$ and $c_2=4$ (two independent vectors in
$E$ determine the third one by multiplication). Thus  $\ci$
is  involutive at $E$ (hence involutive
everywhere, by homogeneity).  Integral manifolds depend on $4$
functions of two variables. 
\end{example}

\section*{Acknowledgements}
These notes are bases on series of lectures
I gave at KIAS    2005 and the seventh international conference
on geometry, integrability and quantization  (respectively March and June 2005). It is a pleasure to
thank my Korean and Bulgarian hosts for these lecture series, especially
 Sung-Eun Koh and Ivalio Mladenov. I also thank T. Ivey, O. Yampolski and
V. Zharnitsky for useful comments.
Supported by NSF grant DMS-0305829


\section*{Bibliography}

 

\begin{thebibliography}{99}
 
\bibitem{BBG} E. Berger, R. Bryant, P. Griffiths,
\emph{  The Gauss equations and rigidity of isometric embeddings},
{  Duke Math. J.} 50 (1983) 803--892.

\bibitem{berger} M. Berger, \emph{  Sur les groupes d'holonomie des variétés riemanniennes}. (French)  C. R. Acad. Sci. Paris  237,  (1953). 472--473


\bibitem{birkhoff} G.  Birkhoff, \emph{Dynamical Systems}, Amer. Math.
Soc., Providence, RI, 1927.

\bibitem{bryantg2} R. Bryant, \emph{  Metrics with exceptional holonomy},
{  Ann. of Math.} 126 (1987), 525--576.
\bibitem{BCG3} R. Bryant, S.-S. Chern, R.B. Gardner, H. Goldschmidt, P.
Griffiths,
\emph{  Exterior Differential Systems}, MSRI Publications, Springer, 1990.



\bibitem{BGY}
R. Bryant, P. Griffiths,D. Yang, 
\emph{  Characteristics and existence of isometric embeddings}.
Duke Math. J. 50 (1983), no. 4, 893--994.

\bibitem{BLZ} Y. Baryishnikov, J. Landsberg, V. Zharnitsky,
\emph{  Periodic orbits in Birkoff billiards},
in preparation.

\bibitem{C1919} E. Cartan, \emph{  Sur les vari\'et\'es de courbure constante d'un
espace euclidien ou non euclidien},
{  Bull. Soc. Math France} 47 (1919) 125--160
and  48 (1920), 132--208; see also pp.~321--432 in
{\it Oeuvres Compl\`etes} Part 3, Gauthier-Villars, 1955.


\bibitem{cartan1931} ---, \emph{  Sur la th\'eorie des syst\`emes en involution et
ses applications \`a la relativit\'e}, {  Bull. Soc. Math. Fr.} 59 (1931),
88--118; see also pp.~1199--1230 in {\it Oeuvres Compl\`etes}, Part 2.


\bibitem{C1945}  ---, \emph{  Les Syst\`emes Ext\'erieurs et leurs
Applications G\'eom\'etriques}, Hermann, 1945.

\bibitem{guill}
V. Guillemin,  
\emph{  The Radon transform on Zoll surfaces}.
Advances in Math. 22 (1976), no. 1, 85--119.

\bibitem{HL} F.  Harvey, H.  Lawson, \emph{  Calibrated geometries},
{  Acta Math.} 148 (1982), 47--157.

\bibitem{IvL}T. Ivey, J.  Landsberg, 
\emph{  Cartan for beginners: differential geometry via moving frames and exterior differential systems}. Graduate Studies in Mathematics, 61. American Mathematical Society, Providence, RI, 2003. xiv+378 pp.



\bibitem{ivrii} V. Ivrii, \emph{  The second term of the spectral asymptotics for the Laplace-Beltrami
 operator on manifolds with boundary},  Functsional. Anal. i Prilozhen., 14, (1980), N 2, 25-34.





\bibitem{joyce} 
M. Gross,  D. Huybrechts, D.  Joyce  \emph{  Calabi-Yau manifolds and related geometries}. Lectures from the Summer School held in Nordfjordeid, June 2001. Universitext. Springer-Verlag, Berlin, 2003. viii+239 pp. 

\bibitem{Ldaewoo} J.  Landsberg
\emph{  Exterior differential systems: A geometric approach to PDE}
in \emph{  Topology and geometry}, proceedings of workshop in pure
math, vol. 17 part III, Korean academic council (1998). pp. 77-100.

\bibitem{rychlik} M.  Rychlik, \emph{  Periodic points of the billiard ball
map in a convex domain}, J. Diff. Geometry {\bf 30} (1989) 191-205.

\bibitem{stojanov} L. Stojanov, \emph{  Note on periodic points of the
billiards}, J. Diff. Geometry {\bf 34} (1991) 835-837.
 

\bibitem{weyl} H. Weyl, \emph{  \"Uber die saymptotische Verteilung der Eigenwerte},
G\"ottinger. Nachr. 110-117 (1911).

\bibitem{mpw} M.  Wojtkovski, \emph{  Two applications of Jacobi fields
to the billiard ball problem}, J. Diff. Geometry {\bf 40} (1994)
155-164.
\end{thebibliography}
\end{document}